\documentclass{amsart}

\usepackage{mathrsfs}

\usepackage{amsfonts,amssymb,amsmath,amsthm}

\usepackage{color}

\usepackage{enumerate}

\usepackage{time}

\def\NN{{\mathbb N}}

\def\RR{{\mathbb R}}

\def\phi{\varphi}

\newcommand{\be}[1]{\begin{equation}\label{#1}}

\newcommand{\ee}{\end{equation}}

\newcommand{\multsum}[2]{\sum_{{\scriptstyle #1}\atop {\scriptstyle #2}}}

\def\cal{\mathscr}

\def\tht{\vartheta}\def\Tht{\Theta}

\def\dd{{\mathrm d}}

\def\NN{{\mathbb N}}

\newtheorem{thm}{Theorem}

\newtheorem{lemma}{Lemma}

\makeatletter

\@addtoreset{equation}{section}
\makeatother

\makeatletter
\def\blfootnote{\xdef\@thefnmark{}\@footnotetext}
\makeatother

\title[Powered numbers]
{On the distribution of powered numbers} 

\author{J. Brüdern and O. Robert}

\begin{document}


\begin{abstract} Asymptotic formulae are established for the number of natural numbers $m$ with 
largest square-free divisor not exceeding $m^\tht$, for any fixed positive parameter $\tht$. Related counting functions are also considered.
\end{abstract}

\maketitle


\section{Introduction}
Motivated by questions about diophantine equations and the $abc$-conjecture, Mazur \cite{M} proposed to smooth out the set of positive $l$-th powers in a multiplicative way, by what he named {\em powered} numbers. To introduce the latter, 
let $k(m)$ denote the largest square-free divisor of the natural number $m$, and let $ i(m) = {\log m}/{\log k(m)}$.
For all $l\in\mathbb N$ one has $i(m^l)\ge l$.  Mazur's powered numbers (relative to $l$) are the numbers $m\in\mathbb N$  with $i(m)\ge l$. Note here  that  $l$ need not be integral in this definition, but for $l\in\mathbb N$ the powered numbers (relative to $l$) contain the $l$-th powers. It is proposed in \cite{M} to replace, within a given diophantine equation, an $l$-th power by the corresponding powered numbers, and to consider the resulting equation between powered numbers as the associated ``rounded'' diophantine equation.  

In this note we analyse the distribution of powered numbers.
We find it is more appropriate to work with the real number $\tht=1/l$. The condition $i(m)\ge l$ is expressed equivalently as $k(m)\le m^\tht$, and for any $\tht>0$ we define the set 
\[
\cal{A}(\tht)=\{m\in \NN\colon  k(m)\le m^{\tht}\}.
\] 
Thus, Mazur's powered numbers (relative to $l$) are exactly the elements of 
$\cal{A}(1/l)$. For analytic approaches to rounded diophantine problems it is indispensable to determine the density of the set $\cal{A}(\tht)$. Our principal goal are asymptotic formulae for the number $S_\tht(x)$ of elements in $\cal{A}(\tht)$ that do not exceed $x$, and for related counting functions.

  It is not difficult to see that for any $0<\tht<1$ the number $S_\tht(x)$  obeys the inequalities
\be{simple} x^\tht \ll S_\tht(x) \ll x^{\tht+\varepsilon} \ee
whenever $\varepsilon$ is a given positive real number and $x$ is large in terms of $\varepsilon$. It now transpires that for $\tht=1/l$ the powered numbers are not much denser than the $l$-th powers. A weaker version of \eqref{simple} occurs in Mazur \cite{M} who
refers to Granville, showing him an ``easy'' argument supposedly confirming the inequalities $x^{\tht-\varepsilon}  \ll S_\tht(x) \ll x^{\tht+\varepsilon}$. In Mazur's article there is no indication how this would go but the simplest argument we know allows one to take $\varepsilon=0$ in the lower bound. To substantiate this claim, fix a number $\tht\in (0,1)$. Define $l\in \mathbb N$ by $1/l<\tht\le 1/(l-1)$, and $t\in\mathbb R$ 
by $l-t=1/\tht$. Then $t\in(0,1]$ and
\be{1/theta}
{1}/{\tht}=t(l-1)+(1-t)l.
\ee
Let $W(x)$ be the number of natural numbers $w$ that have a representation
\be{rep-w}
w=8^{l}n^{l-1}m^{l}
\ee
with $n,m$ square-free and constrained to the intervals
\be{range-n-m}
\frac{1}{2}\,\Big(\frac{x}{8^l}\Big)^{\tht t}<n\le \Big(\frac{x}{8^l}\Big)^{\tht t},\quad \frac{1}{2}\,\Big(\frac{x}{8^l}\Big)^{\tht (1-t)}<m\le \Big(\frac{x}{8^l}\Big)^{\tht (1-t)}.
\ee
The numbers $w$ counted by $W(x)$ uniquely determine the
square-free numbers $n,m$ with \eqref{rep-w} and \eqref{range-n-m}. This follows from unique factorisation. Hence $W(x)$ equals the number of pairs $n,m$ of square-free numbers satisfying \eqref{range-n-m}. In particular, we see that $W(x)\gg x^{\tht}$.

Next, we observe that \eqref{rep-w} gives $k(w)\le 2nm$, and therefore $k(w)\le 2(8^{-l}x)^{\tht}$. By \eqref{1/theta} and \eqref{range-n-m}, we also have $w\ge 2^{-(2l-1)}x$. It follows that $k(w)\le w^{\tht}$. Hence, numbers counted by $W(x)$ are also counted by $S_{\tht}(x)$. This yields $S_{\tht}(x)\ge W(x)$, and the left inequality in \eqref{simple} follows at once.

For the upper bound in \eqref{simple} one may refer to authorities like Tenenbaum \cite[Theorem II.1.15]{T}. However, the shortest argument available to us uses Rankin's trick within the following chain of obvious inequalities:
\[ S_\tht(x) \le \multsum{m\le x}{k(m)\le x^\tht} 1\le \sum_{m=1}^\infty
\Big(\frac{x^\tht}{k(m)}\Big) \Big(\frac{x}{m}\Big)^\varepsilon \ll x^{\tht
+\varepsilon}. \]

A finer analysis of the counting function $S_\tht(x)$ is currently an undesired lacuna in the literature. As a partial remedy we establish an asymptotic formula that should be sufficient for many applications. Let $1\le y\le x$. Following Robert and Tenenbaum \cite{RT2013} we write $N(x,y)$ for the number of $m\in\mathbb N$ with $m\le x$ and $k(m)\le y$. The numbers counted here are their {\em nuclear} numbers. It seems natural to expect that the trivial upper bound
  $S_{\tht}(x)\le N(x,x^{\tht})$ should not be very wasteful.  However, it was conjectured by  Erd\"os (1962)  and  proved by De Bruijn and van Lint (1963) that
\be{erdos62}
\sum_{m\le x}\frac{m}{k(m)}=o\Big( \sum_{m\le x}\frac{x}{k(m)} \Big)\qquad (x\to+\infty).
\ee
For a quantitative version of this estimate, see \cite[Th\'eor\`eme 4.4]{RT2013}. These results suggest that the order of magnitude of $S_{\tht}(x)$ is somewhat smaller than that of $N(x,x^{\tht})$, and this is indeed the case.  

Before we make this precise, 
we recall an estimate for $N(x,y)$. We are primarily interested in the case $y=x^{\tht}$ with $\tht>0$, but we work in the wider range $y>\exp \big((\log x)^{2/3}\big)$. The multiplicative function
\be{psi}
\psi(m)=\prod_{p \mid m}(p+1)\qquad (m\ge 1)
\ee
and the function $F:[0,\infty)\to [0,\infty)$ defined by
\[
F(v)=\frac{6}{\pi^2}\sum_{m\ge 1}\frac{1}{\psi(m)}\min\Big(1,\frac{e^v}{m}\Big)
\]
are featured in the uniform asymptotic formula
\be{nuc2}
N(x,y)=\big(1+o(1)\big)yF\big(\log(x/y)\big)\qquad \big(y> \exp \big((\log x)^{2/3}\big),\thinspace x\to+\infty\big)
\ee
that is contained in \cite[Proposition 10.1]{RT2013}. As we shall see momentarily,
for each pair $\tht, x$ with $0<\tht<1$ and  $x\ge 2$, there is exactly one real number $\alpha=
\alpha_{\tht}(x)>0$ with 
\be{seltsam} 
\sum_{p}\frac{p^{\alpha}\log p}{(p^{\alpha}-1)\big(1+ (p+1)(p^{\alpha}-1)   \big)}
=(1-\tht)\log x. \ee
We are now in a position to state our first result.
\begin{thm}\label{S-theta} Let $0<\tht<1$ fixed. Then, for $x\ge 27$, one has
\be{S7}
S_{\tht}(x)= x^{\tht} F\big( (1-\tht)\log x\big)\frac{  \alpha_{\tht}(x)    }{\tht}   \Bigg(1+   O\Big(  \sqrt{\frac{\log\log x}{\log x}}  \Big)\Bigg) .
\ee
As $x\to \infty$, one also has 
\be{S8}
S_{\tht}(x)=(1+o(1)) x^{\tht} F\big( (1-\tht)\log x\big)\frac{1}{\tht} \Big(\frac{2}{1-\tht}\Big)^{1/2} \big((\log x) \log\log x\big)^{-1/2}.
\ee
\end{thm}

This result calls for several comments. First we take $y=x^\tht$ in \eqref{nuc2} and substitute the resulting equation
\be{Ntheta}
N(x,x^{\tht})=(1+o(1)) x^{\tht} F\big( (1-\tht)\log x\big)\qquad (x\to+\infty)
\ee
within \eqref{S8} to infer that
\[ S_{\tht}(x)=(1+o(1)) N(x,x^{\tht})  \frac{1}{\tht} \Big(\frac{2}{1-\tht}\Big)^{1/2} \big((\log x) \log\log x\big)^{-1/2}      \quad (x\to +\infty).
\]
This is an analogue of \eqref{erdos62}, in quantitative form that is of strength comparable to \cite[Th\'eor\`eme 4.4]{RT2013}.

Our second comment concerns the implicitly defined function $\alpha_\tht(x)$.
It originates in the Dirichlet series
\be{defG}
\cal{G}(s)=  \frac{6}{\pi^2}  \sum_{m\ge 1}\frac{1}{\psi(m)m^s}=\frac{6}{\pi^2}\prod_{p}\left(1+\frac{1}{(p+1)(p^s-1)}\right)
\ee
that converges absolutely in $\text{Re}\, s > 0$, and therefore has no zeros in this half plane. For real numbers $\sigma>0$ we have $\cal G(\sigma)>0$. We may then define
\[
g(\sigma)=\log \cal{G}(\sigma). 
\]
Note that $g$ extends to a holomorphic function on the right half plane, and one computes the logarithmic derivative of $\cal G(s)$ from the
Euler product representation \eqref{defG} to
\[ g'(s) = -  \sum_{p}\frac{p^{s}\log p}{(p^{s}-1)\big(1+ (p+1)(p^{s}-1)   \big)}. \]
Note that the sum on the right hand side here coincides with the sum in \eqref{seltsam}. On differentiating again, it transpires that the real function $g':(0,\infty)\to \mathbb R$ is increasing.  Considering $\sigma\to 0$ and $\sigma\to\infty$ one finds that its range is the open interval $(-\infty,0)$. We conclude that for a given $v>0$ there is a unique positive number  $\sigma_v$ with
\[
g'(\sigma_v)+v=0.
\]
In \cite[Lemme 6.6]{RT2013} it is shown that
\be{estim-sigma-v}
 \sigma_v=(1+o(1))  \sqrt{\frac{2}{v\log v}}\qquad (v\to+\infty).
 \ee
Here we choose $v=(1-\tht)\log x$ and then have 
$\alpha_\tht(x)=\sigma_{v}$. In particular, we see that \eqref{S7} and \eqref{estim-sigma-v} imply \eqref{S8}. Thus, it only remains  to prove \eqref{S7}. 

Our last comment concerns the actual size of $S_\tht(x)$. This requires some more information on the function $F$. 
From \cite[(2.12)]{RT2013}  we have the asymptotic relation
 \be{logF}
 \log F(v)=\big(1+o(1)\big)\sqrt{\frac{8v}{\log v}}\qquad (v\to+\infty) .
 \ee
By inserting \eqref{logF} into \eqref{S8}, we deduce that 
 there exists some positive number  $\beta(x;\tht)>0$ with 
\be{helpme}
S_{\tht}(x)=x^{\tht + \beta(x;\tht)   }\qquad (0<\tht<1,\thinspace x\ge x_0(\tht))
\ee
and the property that for any fixed $\tht\in(0,1)$ one has
\[
\beta(x;\tht)=(1+o(1))\sqrt{\frac{8(1-\tht)}{(\log x)(\log\log x)}}\qquad (x\to+\infty).
\]
Note that \eqref{helpme} yields another proof of \eqref{simple}.

\bigskip
We now turn to local estimates for $S_{\tht}(x)$. In our case, this amounts to comparing the respective behaviour of $S_{\tht}(zx)$ and  $S_{\tht}(x)$ uniformly for large $x$, when $z$ is in some sense sufficiently close to $1$. Such estimates are often obtained with the saddle-point method, and we follow this route here, too. In a suitable range for $z$, the fraction   $S_{\tht}(zx)/S_{\tht}(x)$    may be approximated by a simple function of $z$.

\begin{thm}\label{local-estim}
Let $0<\tht<1$. Then for $x$ large, we have
\[
S_{\tht}(zx)=z^{\tht}S_{\tht}(x)\Bigg(1+   O\Big(  \sqrt{\frac{\log\log x}{\log x}}  \Big)\Bigg) 
\]
uniformly for $z>0$ with $|\log z|\ll \log\log x$.
\end{thm}

Finally, we consider the counting function for a  variation of the powered numbers. For given $\tht\in (0,1)$ and $\Tht\in \RR$, we consider 
$$
S_{\tht,\Tht}(x)=\#\{n\le x\colon k(n)\le n^{\tht}(\log n)^{\Tht}\}.
$$
Note that $S_{\tht}(x)=S_{\tht,0}(x)$. The set of integers such that $k(n)\le n^{\tht}(\log n)^{\Tht}$    plays a prominent role in a forthcoming paper, and therefore, we provide an estimate for $S_{\tht,\Tht}(x)$. It turns out that the conditions $k(n)\le n^{\tht}$ and $k(n)\le n^{\tht}(\log n)^{\Tht}$ are relatively close, and that the ratio $S_{\tht,\Tht}(x)/S_{\tht}(x)$ is roughly of size $(\log x)^{\Tht}$.

\begin{thm}\label{k-log}
Let $0<\tht<1$ and $\Tht\in \RR$ be fixed. Then for $x$ large, one has
\[
S_{\tht,\Tht}(x)=(\log x)^{\Tht}S_{\tht}(x)\Bigg(1+   O\Big(  \sqrt{\frac{\log\log x}{\log x}}    \Big)\Bigg).
\]
\end{thm}

\section{Proof of Theorem 1}

In this section we derive Theorem 1. Before we embark on the main argument, we fix some notation and recall a pivotal result concerned with the distribution of square-free numbers. This involves the function $\psi(m)$ as defined in \eqref{psi}, the M\"obius function $\mu(m)$, and for a parameter $0\le \gamma\le \frac{1}{2}$ at our disposal, the product 
\[
r_\gamma(m)=\prod_{p\mid m}\big(1+4\gamma p^{-1/2}\big).\]
One then has the estimate (\cite[(10.1)]{RT2013})
\be{sum-mu}
\sum_{l\le z} \mu^2(lk)=\frac{6kz}{\pi^2\psi(k)} +O\big( r_\gamma(k)z^{1-\gamma}     \big) \ee
that holds uniformly relative to the square-free number $k$ and the real parameters $z, \gamma$ in the ranges
$z\ge 1$, $0\le \gamma \le \frac12$.

 The first steps of our argument  follow the pattern  laid out in \cite[Sect.\ 10]{RT2013}. Unique factorisation shows that for all natural numbers $n$ there exists exactly one pair of coprime natural numbers $l, m$ with $\mu(l)^2=1$ and $n=lmk(m)$. Note that the two conditions $(l,m)=1$ and $\mu(l)^2=1$ are equivalent to the single condition $\mu(lk(m))^2=1$. Further, one has $k(n)=lk(m)$. With $\tht\in(0,1)$ now fixed, it follows that $S_\tht(x)$ equals the number of $(l,m)\in\mathbb N^2$ satisfying the conditions
\[ \mu(lk(m))^2=1, \quad lmk(m)\le x, \quad lk(m) \le (lmk(m))^\tht. \]
These last three conditions we recast more compactly as
\be{32}
\mu(lk(m))^2=1, \quad lk(m) \le \min (x/m, m^{\tht/(1-\tht)}).
\ee
From now on, the number $\kappa=\tht/(1-\tht)$ features prominently, and we also put $y=x^\tht$. Note that
\be{32a} \min (x/m, m^\kappa) = m^\kappa \;\; \text{ if and only if }\,\, m\le x/y. \ee
Hence, we consider the ranges $m\le x/y$ and $x/y<m\le x$ separately. By \eqref{32} and \eqref{32a}, this leads to the decomposition
\be{33} S_{\tht}(x)=S_{\tht}^{(1)}(x)+S_{\tht}^{(2)}(x) \ee
in which
\[
S_{\tht}^{(1)}(x)  =\!\multsum{m\le x/y}{k(m)\le m^{\kappa}}   \sum_{l\le m^{\kappa}/{k(m)}}\!\mu(lk(m))^2,\quad 
S_{\tht}^{(2)}(x)= \!  \multsum{x/y<m\le x}{mk(m)\le x}   \sum_{l\le {x}/mk(m)}\!\mu(lk(m))^2.
\]
We apply \eqref{sum-mu} with $k=k(m)$ to both inner sums and obtain
\begin{align}
 S_{\tht}^{(1)}(x) &  = \frac{6}{\pi^2}\multsum{m\le x/y}{k(m)\le m^{\kappa}}
 \frac{  m^{\kappa}  }{\psi(m)}+O(R_1), \label{S1} \\
 S_{\tht}^{(2)}(x)  & = \frac{6}{\pi^2}\multsum{x/y<m\le x}{mk(m)\le x} \frac{  x  }{m\psi(m)}+O(R_2) .\label{S2}
\end{align}
where
\[ R_1 = \multsum{m\le x/y}{k(m)\le m^{\kappa}}r_\gamma(m)\Big( \frac{  m^{\kappa  }}{k(m)} \Big)^{1-\gamma}, \qquad
 R_2 = \multsum{x/y<m\le x}{mk(m)\le x} r_\gamma(m)\Big( \frac{ x }{mk(m)} \Big)^{1-\gamma}.\]

It turns out that $R_1$ and $R_2$ are small. In order to couch their estimation, as well as the analysis of other error terms that arise later, under the umbrella of a single treatment, we choose parameters $\gamma$ and $\sigma$ with $0<\gamma<\sigma\le \frac12$ and introduce the series
\be{E} E = \sum_{m=1}^\infty r_\gamma(m)\Big(\frac{y}{k(m)}\Big)^{1-\gamma} \Big(\frac{x/y}{m}\Big)^\sigma. \ee
The conditions $\sigma>\gamma>0$ ensure convergence.  It is routine to show that
\[ R_1+R_2 \ll E. \] 
In fact, by Rankin's trick,
$$ R_1 \le \multsum{m\le x/y}{k(m)\le m^{\kappa}}r_\gamma(m)\Big( \frac{  m^{\kappa  }}{k(m)}\Big)^{1-\gamma} \Big(\frac{x/y}{m}\Big)^\sigma. $$
For $m\le x/y$ one has $m^\kappa \le y$, and it follows that $R_1\le E$.

Likewise, one confirms $R_2\le E$ by observing that $1-\gamma>
\frac12 \ge \sigma$ so that
\be{help} \Big(\frac{x}{mk(m)}\Big)^{1-\gamma} = \Big(\frac{y}{k(m)}\Big)^{1-\gamma}
\Big(\frac{x/y}{m}\Big)^{1-\gamma} \le \Big(\frac{y}{k(m)}\Big)^{1-\gamma}
\Big(\frac{x/y}{m}\Big)^\sigma. \ee

The appearance of $k(m)$ in the summation conditions on the right hand sides of \eqref{S1} and \eqref{S2} is a nuisance, and we proceed by removing these. If the condition $k(m)\le m^\kappa$ is removed from the sum in \eqref{S1} one imports an error no larger than
\[ \multsum{m\le x/y}{k(m)> m^{\kappa}} \frac{  m^{\kappa}  }{\psi(m)} 
 \le \multsum{m\le x/y}{k(m)> m^{\kappa}}
 \frac{  m^{\kappa}  }{k(m)}  
 \le \sum_{m\le x/y} \Big(\frac{m^{\kappa}}{k(m)}\Big)^{1-\gamma} \Big(\frac{x/y}{m}\Big)^\sigma\le E.
 \] 
Here, the last inequality is obtained by the argument that completed the estimation of $R_1$. 
Similarly, if the condition $mk(m)\le x$ is removed from the summation condition in \eqref{S2}, then the resulting error does not exceed
\[ \multsum{x/y<m\le x}{mk(m)> x} \frac{  x  }{m\psi(m)}
\le \multsum{x/y<m\le x}{mk(m)> x} \frac{  x  }{mk(m)}
\le \multsum{m>x/y}{mk(m)> x} \frac{  x  }{mk(m)} = R, \]
say. By Rankin's trick and \eqref{help},
\[ R\le \sum_{m>x/y} \Big(\frac{  x  }{mk(m)}\Big)^{1-\gamma}
\le \sum_{m>x/y} \Big(\frac{  y  }{k(m)}\Big)^{1-\gamma}\Big(\frac{x/y}{m}\Big)^\sigma \le E. \]

Collecting together, we deduce from \eqref{S1} and \eqref{S2} the asymptotic relations
\[
 S_{\tht}^{(1)}(x)   = \frac{6}{\pi^2}\sum_{m\le x/y}
 \frac{  m^{\kappa}  }{\psi(m)}+O(E),  \qquad
 S_{\tht}^{(2)}(x)   = \frac{6}{\pi^2}\sum_{x/y<m\le x} \frac{  x  }{m\psi(m)}+O(E) ,
\]
and by \eqref{33}, we infer that
\[ S_\tht(x) = \frac{6}{\pi^2}\sum_{m\le x} \frac{  m^{\kappa}  }{\psi(m)}\min(1, xm^{-\kappa-1}  ) + O(E). \]
Note that the sum on the right is a partial sum of a convergent series. If one completes the sum, then it is immediate that the error thus imported is bounded by
$R$, and hence by $E$. We have now reached the provisional expansion
\be{S3}S_\tht(x) = \frac{6}{\pi^2}\sum_{m=1}^\infty \frac{  m^{\kappa}  }{\psi(m)}\min(1, xm^{-\kappa-1}  ) + O(E).
\ee

It remains to estimate $E$. In its definition \eqref{E}, we encounter a sum over a multiplicative function, and so
\begin{align*} E =& y^{1-\gamma} (x/y)^\sigma \prod_{p}\Big(1+\frac{p^{\gamma}r(p;\gamma)}{p(p^{\sigma}-1)}\Big)\\ \le& y^{1-\gamma} (x/y)^\sigma\prod_{p}\Big(1+\frac{p^{\gamma}}{p(p^{\sigma}-1)}\Big)\prod_{p}\Big(1+\frac{4\gamma p^{\gamma-\frac{1}{2}}}{p(p^{\sigma}-1)}\Big).
 \end{align*}
Now, since $p^{\gamma-\frac{1}{2}}\le 1$, $\gamma< \sigma$, $p^{\sigma}-1\ge \sigma\log p$, and since the sum $\sum_{p}\frac{1}{p\log p}$ converges, 
 \[
 \prod_{p}\Big(1+\frac{4\gamma p^{\gamma-\frac{1}{2}}}{p(p^{\sigma}-1)}\Big)\le \prod_{p}\Big(1+\frac{4}{p\log p}\Big)\ll 1.
 \]
As on an earlier occasion, we write $v= (1-\tht)\log x = \log \frac{x}{y}$, and then choose $\sigma=\sigma_v$ and $\gamma=\sigma_v-\frac{1}{\log y}$. By \eqref{defG} one has 
 \[
 \prod_{p}\Big(1+\frac{p^{\gamma}}{p(p^{\sigma}-1)}\Big)=\frac{\pi^2}{6}\cal{G}(\sigma)\prod_p\Big(1+\frac{(p^{\gamma}-1)(p+1)+1}{p+p(p+1)(p^{\sigma}-1)}    \Big).
 \]
 Now, on recalling \eqref{estim-sigma-v}, 
 \[
 \prod_{p\le e^{1/\sigma}}\Big(1+\frac{(p^{\gamma}-1)(p+1)+1}{p+p(p+1)(p^{\sigma}-1)}    \Big)\le \prod_{p\le e^{1/\sigma}}\Big(1+\frac{1}{p}\Big)\ll \frac{1}{\sigma}\ll (v\log v)^{1/2}
 \]
while one also has
 \[
 \prod_{p> e^{1/\sigma}}\Big(1+\frac{(p^{\gamma}-1)(p+1)+1}{p+p(p+1)(p^{\sigma}-1)}    \Big)\le \prod_{p> e^{1/\sigma}}\left(1+\frac{1}{p^{1+\sigma-\gamma}}\right)\ll \zeta\Big(1+\frac{1}{\log y}\Big)\ll \log y.
 \]
On collecting together, this shows that
 \[
  E\ll y^{1-\sigma_v}e^{ v\sigma_v}\cal{G}(\sigma_v)(v\log v)^{1/2}\log y.
\]
From \cite[(2.11)]{RT2013} we deduce that
 \[
 F(v)\asymp   \Big(\frac{\log v}{v}\Big)^{1/4} e^{v\sigma_v} \cal{G}(\sigma_v)    \qquad (v\ge 2),
 \]
and hence   
\[  
E \ll y^{1-\sigma_v}F(v) v^{3/4} (\log v)^{1/4} \log y.
 \] 
 With the choice of $y$ and $v$, one has $\sigma_v=\alpha_{\tht}(x)$.  Moreover,  $\log y$ and $v$ have the order of magnitude $\log x$ so that the last inequality now reads
 \be{Enew}
 E\ll x^{\tht} F((1-\tht)\log x)x^{-\tht \alpha_{\tht}(x)}(\log x)^{7/4}(\log\log x)^{1/4}.
 \ee

\medskip
Our final task is to compare our estimate for $E$ with the size of the sum on the right of \eqref{S3}.

Recall that in view of  \eqref{simple} and \eqref{Ntheta},   $S_{\tht}(x)$ and $N(x,x^{\tht})$ are of comparable size. In
order to mimick the estimate \eqref{Ntheta},
we introduce the function
\[
H_{\tht}(x)=\frac{6}{\pi^2x^{\tht}}\sum_{m=1}^\infty \frac{  m^{\kappa}  }{\psi(m)}\min(1, xm^{-\kappa-1}  )
\]
so that \eqref{S3} now reads
\[
S_\tht(x) =x^{\tht}H_{\tht}(x)+O(E).
\]

Our aim is to give an estimate of $H_\tht(x)$ by using the saddle-point method, and to describe more precisely  $ H_\tht(x) /F((1-\tht)\log x)$ as $x\to+\infty$.

\begin{lemma}\label{lemmaH}
Fix $\tht$ with $0<\tht<1$. Then for $x\ge 27$ one has 
\[
H_\tht(x)=F((1-\tht)\log x)\frac{\alpha_{\tht}(x)}{\tht}\Bigg(1+O\Big(\sqrt{\frac{\log\log x}{\log x}}\Big)\Bigg).
\]
\end{lemma} \noindent
\textit{Proof.} The argument is modelled on  \cite[Section 8]{RT2013}. Recall the identity 
\[
\frac{1}{2\pi \mathrm{i}}\int_{\sigma+\mathrm{i}\RR}\frac{y^{s}}{s(1-s)}\,\dd s=\min(1,y),\qquad (0<\sigma<1,\thinspace y>0).
\]
We then have the integral representation
\begin{align*}
H_\tht(x)&=\frac{1}{x^{\tht}}\frac{6}{\pi^2}\sum_{m\ge 1} \frac{  m^{\kappa}  }{\psi(m)}\frac{1}{2\pi \mathrm{i}}\int_{\sigma+\mathrm{i}\RR}\frac{1}{s(1-s)}  \left(\frac{x}{m^{1/(1-\tht)}}\right)^s\dd s
\\
&=  \frac{1}{2\pi \mathrm{i}}\int_{\sigma+\mathrm{i}\RR}  \cal{G}\Big(\frac{s-\tht}{1-\tht}\Big)x^{s-\tht}     \frac{\dd s}{s(1-s)} 
\\
&=  \frac{1}{2\pi }\int_{\RR}  \cal{G}\Big(\frac{\sigma+\mathrm{i}t -\tht}{1-\tht}\Big)x^{\sigma+\mathrm{i}t-\tht}     \frac{\dd t}{(\sigma+\mathrm{i}t)(1-\sigma-\mathrm{i}t)} .
\end{align*}
After a linear change of variable in $t$, we arrive at
\[H_\tht(x)=
  \frac{1}{2\pi }\int_{\RR}  \cal{G}\Big(\frac{\sigma-\tht}{1-\tht}+\mathrm{i}t \Big)x^{\sigma+\mathrm{i}(1-\tht)t-\tht}     \frac{(1-\tht)\,\dd t}{(\sigma+\mathrm{i}(1-\tht)t)(1-\sigma-\mathrm{i}(1-\tht)t)} .
\]

Recall again that $v=(1-\tht)\log x$, and that $\alpha_{\tht}(x)=\sigma_v$. We take $\sigma=\tht+(1-\tht)\sigma_v$. For large $x$ one then has $0<\sigma<1$, and  the previous formula for $H_\tht(x)$ becomes
\[
H_\tht(x)= \frac{1}{2\pi }\int_{\RR} \frac{ \cal{G}(\sigma_v+\mathrm{i}t) e^{(\sigma_v+\mathrm{i}t)v} }
{\big(\tht+(1-\tht)(\sigma_v+\mathrm{i}t)\big)\big( 1-\sigma_v-\mathrm{i}t      \big) } \, {\dd t}.
\]
After truncation, we have 
\[
H_\tht(x)=\frac{1}{2\pi } \int_{-v^2}^{v^2}
\frac{ \cal{G}(\sigma_v+\mathrm{i}t) e^{(\sigma_v+\mathrm{i}t)v} }
{\big(\tht+(1-\tht)(\sigma_v+\mathrm{i}t)\big)\big( 1-\sigma_v-\mathrm{i}t      \big) } \, \dd t
+ 
O\left(   \frac{\cal{G}(\sigma_v) e^{v\sigma_v} }{v^2}   \right).
\]
Moreover, following the proof of \cite[Th\'eor\`eme 8.6]{RT2013}, we set
\[
\eta_v=(\log v)/\sqrt{g''(\sigma_v)}\asymp (\log v)^{3/4}/v^{3/4}
\]
and recall \cite[Lemme 8.5]{RT2013}, asserting that for some $c>0$ we have
\[
|\cal{G}(\sigma_v+\mathrm{i}t)|\ll \cal{G}(\sigma_v)e^{-c (\log v)^2}\qquad (\eta_v\le  |t|\le \exp((\log v)^{38/37}).
\]
It now follows that
\[
H_\tht(x)=\frac{1}{2\pi } 
\int_{-\eta_v}^{\eta_v}
\frac{ \cal{G}(\sigma_v+\mathrm{i}t) e^{(\sigma_v+\mathrm{i}t)v} }
{\big(\tht+(1-\tht)(\sigma_v+\mathrm{i}t)\big)\big( 1-\sigma_v-\mathrm{i}t      \big) } \, {\dd t}
+ 
O\left(   \frac{\cal{G}(\sigma_v) e^{v\sigma_v} }{v^2}   \right).
\]
Setting
\[
 D_{m}=v^{(m+1)/2}(\log v)^{(m-1)/2}\qquad (m\ge 0),
\]
we have
\[
H_\tht(x)=\frac{ \cal{G}(\sigma_v) e^{v\sigma_v}       }{2\pi }  \int_{-\eta_v}^{\eta_v}  \Upsilon (t)e^{-g''(\sigma_v)t^2/2} \dd t+   O\left(   \frac{\cal{G}(\sigma_v) e^{v\sigma_v} }{v^2}   \right)
\]
where
\[
\Upsilon(t)= \frac{e^{-\mathrm{i}g'''(\sigma_v)t^3/6+O(t^4D_4)}}{  \big(\tht+(1-\tht)(\sigma_v+\mathrm{i}t)\big)\big( 1-\sigma_v-\mathrm{i}t      \big) }.   
\]
By Taylor expansion, we infer
\[ \Upsilon(t) = \frac{\mathrm Z(t)}{(\tht+(1-\tht)\sigma_v)(1-\sigma_v) }\]
where
\[ \mathrm Z(t)=
1+\mathrm{i}t \Big(\frac1{1-\sigma_v}-\frac{1-\tht}{\tht+(1-\tht)\sigma_v}     \Big)      -{\mathrm{i}}g'''(\sigma_v)\frac{t^3}6         +O(t^2+D_4t^4+D_3^2t^6)  .
\]
Now, still following the pattern of  the proof of \cite[Th\'eor\`eme 8.6]{RT2013}, one is lead to 
\[
H_\tht(x)=\frac{x^{(1-\tht) \alpha_{\tht}(x)  } \cal{G}\big(\alpha_{\tht}(x)\big)    }{\tht \sqrt{2\pi g''\big( \alpha_{\tht}(x)   \big)    }}\Bigg(1+O\Big(\sqrt{\frac{\log\log x}{\log x}}\Big)\Bigg).
\]
We omit the details. From \cite[Th\'eor\`eme 8.6]{RT2013} we import the relation
 \[
 F(v)=\frac{e^{v\sigma_v } \cal{G}(\sigma_v)   }{\sigma_v\sqrt{2\pi g''(\sigma_v)}}\left(1+O\left(    \sqrt{\frac{\log v}{v}}          \right)\right)\qquad (v\ge 2),
 \]
and the lemma follows.
\qed

\smallskip
We may now complete the proof of Theorem 1.   Using the  lemma and \eqref{estim-sigma-v}, we obtain
\[
H_{\tht}(x)\asymp F((1-\tht)\log x) \big((\log x)\log\log x\big)^{-1/2},
\]
so that the estimate \eqref{Enew} implies
\[
E\ll     x^{\tht}H_{\tht}(x)x^{-\tht\alpha_{\tht}(x)}(\log x)^{9/4}(\log\log x)^{3/4}.
\]
We then have
\[
S_{\tht}(x)=x^{\tht}H_{\tht}(x)\big(1+O\big(    x^{-\tht\alpha_{\tht}(x)}(\log x)^{9/4}(\log\log x)^{3/4}           \big)\big).
\]
We may  now  replace $H_{\tht}(x)$ by the estimate from Lemma  \ref{lemmaH}. Since 
\[
x^{-\tht\alpha_{\tht}(x)}(\log x)^{9/4}(\log\log x)^{3/4}\ll \sqrt{\frac{\log\log x}{\log x}},
\]
this yields \eqref{S7}. As remarked earlier, \eqref{S8} follows from \eqref{estim-sigma-v} and \eqref{S7}.  The proof of Theorem 1 is complete.

\section{Proof of Theorem \ref{local-estim} }

Subject to the hypotheses of Theorem \ref{local-estim}, when $x$ is large, one has $\log zx\asymp \log x$. Hence, Theorem 1 implies that
\be{Szx}
S_{\tht}(zx)=(zx)^{\tht} F\big( (1-\tht)\log zx\big)\frac{  \alpha_{\tht}(zx)    }{\tht}   \Bigg(1+   O\Big(  \sqrt{\frac{\log\log x}{\log x}}  \Big)\Bigg) 
\ee
holds uniformly for $|\log z|\ll \log\log x$. We recall \cite[Proposition 8.7]{RT2013}. This asserts that uniformly for $|t|\le v^{3/4}(\log v)^{1/4}$, one has
$$
F(v+t)=e^{t\sigma_{v}}F(v)\Big(1+O\Big(\frac{\log v + t^2/v}{\sqrt{v\log v}}\Big)   \Big).
$$
 Using this estimate with $v=(1-\tht)\log x$ and $t=(1-\tht)\log z$, one finds
\[
F\big( (1-\tht)\log zx\big)=F\big( (1-\tht)\log x\big)z^{(1-\tht)\alpha_{\tht}(x)} \Bigg(1+   O\Big(  \sqrt{\frac{\log\log x}{\log x}}  \Big)\Bigg).
\]
Moreover, in the ranges for $x$ and $z$ considered here, one has
\[
z^{(1-\tht)\alpha_{\tht}(x)}=1+O\big( \alpha_{\tht}(x)\log z  \big),
\]
which yields an admissible error term.

Finally, \cite[(8.9)]{RT2013} implies that uniformly for $|t|\le v/2$ one has
$$
\sigma_{v+t}=\sigma_v\big(1+ O(|t|/v)\big).
$$
 Hence  
\[
\alpha_{\tht}(zx) =\alpha_{\tht}(x) \Big(1+   O\Big( \frac{\log z}{\log x}  \Big)\Big)
\]
 which again leads to an admissible error term. Inserting these estimates in \eqref{Szx} completes the proof of Theorem \ref{local-estim}.

\medskip 
It may be worth pointing out that the above argument actually proves a little more. A close inspection of the proof of Theorem \ref{local-estim} shows that the   estimate
 \[
S_{\tht}(zx)=z^{\tht+(1-\tht)\alpha_0(x)}S_{\tht}(x)\Bigg(1+   O\Big(  \sqrt{\frac{\log\log x}{\log x}} +\frac{(\log z)^2}{(\log x)^{3/2}(\log\log x)^{1/2}} \Big)\Bigg) 
\]
holds uniformly in the range $z>0$, $x>27$,  $|\log z|\le (\log x)^{3/4}(\log\log x)^{1/4}$.

\section{Proof of  Theorem \ref{k-log}}
 Before proving Theorem \ref{k-log}, we briefly sketch the main steps. We first choose a suitable real number $U=U(x)$ such that $\log U=(\log x)(1+o(1))$, and count the integers $n$ not exceeding $x$ sucht that $k(n)\le n^{\tht}(\log U)^{\Tht}$. The first step is to show that the number of these integers is essentially $S_{\tht}(x)$ multiplied by $(\log x)^{\Tht}$. The second step is to prove that the number of these integers is close to $S_{\tht,\Tht}(x)$.
 
 In light of this description, for any $x\ge 1$ and any $z>0$, we set
 $$
 B(x,z)=\#\{n\le x\colon k(n)\le z n^{\tht}\}.
 $$

\begin{thm}\label{logx}Let $0<\tht<1$ be fixed. Then for $x$ large, one has
\[
B(x,z)=z S_{\tht}(x)\Bigg(1+   O\Big(  \sqrt{\frac{\log\log x}{\log x}}  \Big)\Bigg) 
\]
uniformly for $z>0$ with $|\log z|\ll \log\log x$.
\end{thm}

\textit{Proof. }We follow very closely the proof of Theorem 1 which corresponds to the case $z=1$. 
We redefine the meaning of $y$, now set to $y=x^{\tht}z$, and keep the notation $\kappa$ and $E$.  Note that in the sequel of this proof, the error term $E$ is to be interpreted with the current specific choices for $x$ and $y$.

For any  $n\ge 1$,  recall that there is a unique pair $l,m$ with $n=lmk(m)$ and $\mu^2\big(lk(m)\big)=1$. The conditions $n\le x$ and $k(n)\le z^{\tht}$ become
\[
lmk(n)\le x\quad\mbox{ and }\quad lk(m)\le m^{\kappa}z^{1+\kappa}.
\]
By keeping the same dichotomy $m\le x/y$ and $m>x/y$ it is straightforward that
\[
B(x,z)=\frac{6}{\pi^2}\sum_{m=1}^{\infty}\frac{1}{\psi(m)}\min\big(m^{\kappa}z^{1+\kappa},x/m   \big)+O(E).
\]
Moreover, \eqref{S3} with the parameter $xz^{-\kappa-1}$ reads
\[
S_{\tht}(xz^{-\kappa-1})=\frac{6}{\pi^2}\sum_{m=1}^{\infty}\frac{1}{\psi(m)}\min\big(m^{\kappa},xz^{-\kappa-1}/m   \big)+O(z^{-(1+\kappa)(1-\gamma)}E).
\] 
Hence, one has
\[
B(x,z)=z^{1+\kappa}S_{\tht}(xz^{-\kappa-1})+O(z^{(1+\kappa)\gamma}E).
\]

Still choosing $\gamma=\sigma_v-\frac{1}{\log y}$ and $\sigma=\sigma_v$, we have $z^{(1+\kappa)\gamma}\ll 1$. Moreover, \eqref{Enew}  implies that for some $c>0$ one has
\[
E\ll x^{\tht} F\big( (1-\tht)\log x \big)\exp\Big(-c\sqrt{\frac{\log x}{\log\log x}}\Big).
\]
This is sufficient to ensure that 
\[
E\ll  z S_{\tht}(x)    \sqrt{\frac{\log\log x}{\log x}}.
\]
Finally we estimate the term $S_{\tht}(xz^{-\kappa-1})$ by Theorem \ref{local-estim}. Inserting this in the estimate for $B(x,z)$ and noticing that in the main term the exponent  $1+\kappa-(1+\kappa)\tht $ of $z$  is equal to $1$ gives the expected result. \qed
\smallskip

We may now complete the proof of Theorem \ref{k-log}. It is sufficient to prove the result for $\Tht\neq 0$, since in the case $\Tht=0$ one has $S_{\tht,\Tht}(x)=S_{\tht}(x)$. We put 
\[C=1+(1+|\Tht|)/\tht \quad\text{and}\quad
U=U(x)=x(\log x)^{-C}.
\]
Note that we have $C>0$, $C\tht \ge 1$ and $\Tht+ C\tht \ge 1$.

\smallskip
First consider the case $\Tht>0$. Any integer $n$ counted by $S_{\tht,\Tht}(x)$ satisfies $k(n)\le n^{\tht}(\log x)^{\Tht}$, whence $S_{\tht,\Tht}(x)\le B\big(x,(\log x)^{\Tht}\big)$. Now, a lower bound is obtained by noticing that the set of integers $U<n\le x$ counted in $S_{\tht,\Tht}(x)$ contains the integers such that $k(n)\le n^{\tht}(\log U)^{\Tht}$. These deliberations yield the inequalities
\[
  B\big(x,(\log U)^{\Tht}\big)-B\big(U,(\log U)^{\Tht}\big)            \le   S_{\tht,\Tht}(x)\le B\big(x,(\log x)^{\Tht}\big).                                    
\]
Now, using Theorem \ref{logx} to estimate $ B\big(x,(\log U)^{\Tht}\big)$ and $ B\big(x,(\log x)^{\Tht}\big)$, and then replacing $(\log U)^{\Tht}$ by $(\log x)^{\Tht}$ at the price of an admissible error term, one obtains the main term to estimate $ S_{\tht,\Tht}(x)$.  For the remaining term,  Theorem \ref{logx} and the definition of $U$  imply that
\[
B\big(U,(\log U)^{\Tht}\big)\ll (\log U)^{\Tht}S_{\tht}(U)\ll (\log x)^{\Tht}S_{\tht}(U).
\]
Finally, Theorem \ref{local-estim} with the choice $z=(\log x)^{-C}$ implies that
\[
S_{\tht}(U)\ll z^{\tht}S_{\tht}(x)\ll      (\log x)^{-\tht C}S_{\tht}(x)       \ll (\log x)^{-1}S_{\tht}(x).
\]
Gathering these estimates, we obtain an upper bound for $B\big(U,(\log U)^{\Tht}\big)$ that provides an admissible error term, and therefore proves the result for $\Tht>0$.

\smallskip
The case $\Tht<0$ is very similar. Any integer counted by $S_{\tht,\Tht}(x)$ either satisfies $n\le U$ with $k(n)\le n^{\tht}$, or $U<n\le x$ with $k(n)\le n^{\tht}(\log U)^{\Tht}$. Moreover, any integer $n\le x$ such that $k(n)\le n^{\tht}(\log x)^{\Tht}$ is counted in $S_{\tht,\Tht}(x)$. We conclude that 
\[
 B\big(x,(\log x)^{\Tht}\big)   \le   S_{\tht,\Tht}(x)\le           B\big(x,(\log U)^{\Tht}\big)+ S_{\tht}(U).
\]
As before, one uses Theorem \ref{logx} to estimate  $B\big(x,(\log U)^{\Tht}\big)$ and   $B\big(x,(\log x)^{\Tht}\big)$. This provides the main term. Finally, Theorem \ref{local-estim}  with $z=(\log x)^{-C}$ provides 
\[
S_{\tht}(U)\ll z^{\tht}S_{\tht}(x)\ll      (\log x)^{-\tht C}S_{\tht}(x)       \ll (\log x)^{\Tht -1}S_{\tht}(x),
\]
which yields an admissible error term.

\vspace{2ex}\noindent
 J\"org Br\"udern \\
 Universit\"at G\"ottingen\\
Mathematisches Institut\\
Bunsenstrasse 3--5\\
D 37073 G\"ottingen \\
Germany\\
jbruede@gwdg.de\\
[2ex]
Olivier Robert\\
 Universit\'e de Lyon,  Universit\'e de Saint-\'Etienne\\
Institut Camille Jordan CNRS UMR 5208\\
25, rue Dr R\'emy Annino\\
F-42000, Saint-\'Etienne\\
 France\\
olivier.robert@univ-st-etienne.fr

\end{document}